\documentstyle{amsppt}
\loadbold
\TagsOnRight
\hsize 30pc
\vsize 47pc
\magnification 1100

\redefine\w{\wedge}
\redefine\g{\frak g}

\redefine\a{\frak a}
\redefine\s{\frak s}

\redefine\C{{\Bbb C}}

\redefine\H{{\Cal H}}

\redefine\PC{\Bbb{P}_n(\Bbb C)}

\redefine\PR{\Bbb{P}_n(\Bbb R)}

\redefine\PO{\Bbb{P}_2(\Bbb O)}

\redefine\g{\frak g}
\redefine\k{\frak k}
\redefine\h{\frak h}

\redefine\m{\frak m}

\redefine\u{\frak u}
\redefine\t{\frak t}
\redefine\su{{\frak{su}}}
\redefine\spin{{\frak{spin}}}

\redefine\Ad{{\text{Ad}}}
\redefine\so{{\frak{so}}}
\redefine\sp{{\frak{sp}}}

\redefine\C{\Bbb C}
\redefine\R{\Bbb R}

\redefine\t{\frak t}
\redefine\a{\frak a}
\redefine\Spin{{\text{Spin}}}
\redefine\Sp{{\text{Sp}}}
\redefine\O{{\text{O}}}
\redefine\U{{\text{U}}}
\redefine\T{{\text{T}}}
\redefine\E{{\text{E}}}

\redefine\SO{{\text{SO}}}
\redefine\SU{{\text{SU}}}
\redefine\G2{{\text{G}}_2}
\redefine\F4{{\text{F}}_4}

\redefine\Si{{\Sigma}}
\topmatter
\title Polar and coisotropic actions on K\"ahler manifolds \endtitle
\author Fabio Podest\`a and Gudlaugur Thorbergsson\endauthor
\address \endaddress
\newcount\minutes
\newcount\scratch

\thanks Part of the work on this paper was done during a visit of the second
author at the University of Florence that was financially supported
by G.N.S.A.G.A. - I.N.d.A.M.
\endthanks
 \abstract The main result of the paper is that a polar action on a
compact irreducible homogeneous K\"ahler manifold is coiso\-tropic.
This is then used to give new examples of polar
actions and to classify
coisotropic and polar actions on quadrics.
\endabstract

\leftheadtext{\smc }
\rightheadtext{\smc }

\endtopmatter

\document

\noindent{\bf \S 1. Introduction}
\bigskip
The aim of the present paper is to investigate the relationship
between polar and coisotropic actions on compact K\"ahler manifolds.\par
The action of a compact Lie group $G$ of isometries on a Riemannian manifold
($M,g$) is called {\it polar\/} if there
exists a properly embedded submanifold $\Sigma$ which meets every
$G$-orbit and is orthogonal to the $G$-orbits in all common
points. Such a submanifold $\Sigma$
is called a {\it section\/} (see [PT1,2]) and if it is
flat, the action is
called {\it hyperpolar\/}. It is of course meaningful to relax the
definition of
polar actions and not require that the section is properly embedded. Our
results are correct
for weakly polar actions in this sense, making Theorems 1.1 and 1.3
stronger, but Theorem 1.2 weaker.\par
If ($M,g$) is a compact K\"ahler manifold with K\"ahler form $\omega$ and
$G$
is a compact subgroup of its full isometry group,
then the
$G$-action is called {\it multiplicity-free\/} ([GS2]) or {\it
coisotropic\/}
([HW]) if the principal $G$-orbits are
coisotropic with respect to $\omega$. Notice that the existence of one
coisotropic principal $G$-orbit implies the same property for all
principal $G$-orbits, see [HW].  \par

In this paper we shall consider the special case of $M$ being a simply
connected, compact homogeneous K\"ahler manifold; such manifolds are
diffeomorphic to quotient spaces $L/K$, where $L$
is a compact semisimple Lie group and $K$ is
the centralizer in $L$ of some torus. We also recall that
a compact homogeneous K\"ahler manifold cannot be written as a nontrivial
Riemannian product of two Riemannian manifolds if and only if its full
isometry group is a compact simple Lie group;
in this last case the manifold is
simply connected (see e.g. [On], p.~238 ff).
\par
In section two we will prove the following
\medskip
\proclaim{Theorem 1.1} A polar action on an irreducible compact
homogeneous K\"ahler manifold is coiso\-tropic.\endproclaim
\medskip
 It is obvious that the theorem is not true without
the assumption that $M$ is irreducible. The converse of Theorem 1.1 is not
true. In fact, we will give an
example of a coisotropic action
on $\PC$ which is not polar in section three
and there are further such examples
on quadrics in Theorem 1.3.
In certain situations, we can prove that
coisotropic actions are polar and
we will use this to give many new examples of polar actions. We will for
example show that every nonoriented
compact surface $N$ is the section  of a polar action on some four
dimensional compact
K\"ahler manifold. In section three we will prove the following
\medskip
\proclaim{Theorem 1.2} An effective and isometric action of the (real) torus
$T^n$
on a 
compact  K\"ahler manifold $M$ with positive Euler characteristic and
complex dimension
$n$ is polar.\endproclaim
\medskip

The condition on the Euler characteristic is to ensure that the action is
Poisson; see the proof in section three. We also remark that this condition
is not restrictive, since a compact, symplectic $2n$-dimensional
manifold which is acted on effectively by an $n$-dimensional torus in a
Poisson
fashion has positive Euler characteristic (see [De]).

\medskip

As an application of Theorem 1.1, we will classify polar actions on quadrics
in section four. We do this by first classifying the coisotropic actions and
then 
deciding which of them are polar. An interesting consequence of this
classification 
is that polar actions on quadrics are hyperpolar. This is in contrast with
complex 
projective spaces or more generally to rank one symmetric spaces that admit
many 
polar actions that are not hyperpolar (see [PTh]). Another interesting
consequence of the
classification is that all coisotropic actions on quadrics are related to
transitive actions on spheres. The classification of coisotropic and polar
actions on quadrics is given in the following

\medskip

\noindent{\bf Theorem 1.3}. {\it Let $G$ be a connected compact
Lie subgroup of} $\SO(n)$ {\it ($n \geq 5$) and let $Q_{n-2}$ be the
quadric} 
$\SO(n)/\SO(2)\times \SO(n-2)$. {\it Then the action of $G$ on $Q_{n-2}$ is
nontransitive and coisotropic if and only if one of the following occurs.}
\roster
\item "(i)"  {\it The action of $G$ on $\R^n$ is irreducible and $G$ is
conjugate to}
$$\SU({n\over 2}),\;\U({n\over 2}),\; \Sp(1)\cdot \Sp({n\over 4}),\;
\T^1\cdot \Sp(n/4),\; \Spin(9).$$
{\it Only the action of the last one of these groups is not polar on
$Q_{n-2}$. 
\item "(ii)" The action of $G$ on $\R^n$ has a one dimensional fixed point
set and the group $G$ as a subgroup of} $\SO(n-1)$  {\it is conjugate to one
of the following }
$$\eqalign{{}& \SO(n-1),\;\G2,\;\Spin(7),\cr
{}&\SU({{(n-1)}\over 2}),\;\U({{(n-1)}\over 2}),\;
\Sp(1)\cdot\Sp({{(n-1)}\over 4}),\cr
{}& T^1\cdot\Sp({{(n-1)}\over 4}),\;\Spin(9).}$$ {\it
Only the actions of the first three of these groups are polar on $Q_{n-2}$.
\item "(iii)" The action of $G$ on $\R^n$ is reducible and $\R^n$ is the sum
of two irreducible nontrivial
$G$-submodules $V_1,V_2$; the group $G$ splits as the
product $G_1\times G_2$ with} $G_i \subset \SO(V_i)$, $i = 1,2$ {\it and
each
$G_i$ is conjugate in} $\SO(V_i)$ {\it to one of the following}
$$\SO(p),\; \G2,\; \Spin(7),\; \Sp(1)\cdot\Sp(p),\; \Spin(9)$$
{\it for suitable $p \in {\Bbb N}$. Only the actions
obtained choosing each $G_i$, $i = 1,2$, among the first three of
these groups are
polar on $Q_{n-2}$.}
\endroster
{\it Moreover if the $G$-action on $Q_{n-2}$ is polar, then it is
hyperpolar.}
\bigskip
We now recall some basic definition and
results from the paper [HW]. A submanifold $W$
of a symplectic manifold $(M,\omega)$ is called {\it coisotropic\/} if
$$T_pW^{\perp_\omega}\subseteq T_pW $$
for all $p$ in $W$. Here $T_pW^{\perp_\omega}$ denotes the subspace of
$T_pM$ that
is 
$\omega$-orthogonal to $T_pW$. If $M$ is a K\"ahler
manifold and
$\omega$ its K\"ahler form,
then it is easy to see that a submanifold $W$ of $M$ is coisotropic if and
only if 
$$J(N_pW)\subseteq T_pW$$
for all $p$ in $W$, where $J$ denotes the complex structure of $M$, and
$N_pW$ the normal space of
$W$ in $p$. One calls a real
subspace $V$ in a complex vector space with
a Hermitian scalar
product {\it totally real\/}
if $V$ is perpendicular to $J(V)$. It follows from this discussion that a
submanifold $W$ in a K\"ahler manifold
is coisotropic if and only if all its normal spaces are totally real.

Let $(M,\omega)$ be a symplectic manifold. Then one can associate to every
function $f$ in $C^\infty(M)$
its symplectic gradient $X_f$ defined as the unique vector field on $M$ such
that 
$$df(Y)=\omega(X_f,Y)$$
for every vector field $Y$ on $M$. If $f$ and $g$ are two functions in
$C^\infty$, then one defines their {\it
Poisson bracket} as $\{f,g\}=\omega(X_f,X_g)$. It follows that $C^\infty(M)$
with the Poisson bracket is a Lie
algebra.

Now let $G$ act on $M$ preserving the symplectic structure $\omega$. Then
the action is called {\it Poisson}
if there is a Lie algebra homomorphism $\lambda:{\frak g}\to C^\infty(M)$
such that $X_{\lambda(\xi)}$ agrees
with the infinitesimal action of $\xi$ on $M$. The {\it moment map} of a
Poisson action is defined as
$$\Phi:M\to \frak g^*;\quad \Phi(p)(\xi)=\lambda(\xi)(p).$$
We have $\Phi(gp)=\Ad^*(g)\Phi(p)$, i.e., the moment map is
$G$-equivariant.

Let $G$ be a compact group acting isometrically on a
compact K\"ahler manifold. This
action is automatically holomorphic by a theorem of Kostant
(see [KN], vol.~I, p.~247)
and it induces by compactness of $M$
an action of the complexified group $G^{\Bbb C}$ on $M$. We say that $M$ is
{\it $G^{\Bbb C}$-almost homogeneous}
if $G^{\Bbb C}$ has an open orbit in $M$. If all Borel subgroups of $G^{\Bbb
C}$
act with an open  orbit on $M$, then the open orbit $\Omega$ in $M$ under
the action of $G^{\Bbb C}$ is called a
{\it spherical homogeneous space} and $M$ is called a {\it spherical
embedding of $\Omega$}. In general, given a complex
Lie group $U$ and a closed subgroup $H \subset U$, the pair ($U,H$) is
called a {\it spherical pair\/} if every Borel subgroup of $U$ has an
open orbit in $U/H$; all such pairs ($U,H$), where $U$ is semisimple and $H$
is reductive, have been classified by Brion ([Br2]).

One of the basic results on spherical embeddings is that if the orbit
$\Omega$ of $G^{\Bbb C}$ is a spherical homogeneous
space, then $G^{\Bbb C}$ and all Borel subgroups of $G^{\Bbb C}$ have
finitely many orbits in $M$, see [Br1].

The following theorem  is proved in [HW], p.~275.

\proclaim{Theorem 1.4 (Equivalence Theorem of [HW])} Let $M$ be a connected
compact K\"ahler manifold with an isometric action of a connected compact
group $G$ that is also Poisson. Then the following conditions
are equivalent:
\roster
\item"(i)" The space $C^\infty(M)^G$ of $G$-invariant functions on $M$  is
abelian with respect to the Poisson bracket.
\item" (ii)" The $G$-action is coisotropic.
\item"(iii)" The
cohomogeneity of the $G$-action is equal to the difference between the
rank of $G$ and the rank of a regular isotropy subgroup of $G$.
\item"(iv)" The moment map $\Phi:M\to {\frak g}^*$ separates orbits.
\item"(v)" The K\"ahler manifold $M$ is projective algebraic, $G^{\Bbb
C}$-almost homogeneous and a spherical embedding of
the open $G^{\Bbb C}$-orbit.
\endroster
\endproclaim
\medskip
We remark here that the conditions (i) to (iii) are equivalent even without
the hypothesis of compactness for $M$ (see [HW]).\par
Under the hypothesis of Theorem 1.4, the following convexity theorem holds
(see [GS1] and [Kir]): if $T$ denotes a maximal torus of $G$,
$\t^*$ is the fixed point set of $T$ in $\g^*$ and $\t^*_+$ denotes
a fundamental domain for the action of the Weyl group on $\t^*$, then
the image $\Phi(M)$ under the moment mapping $\Phi$ intersects $\t^*_+$
in a convex polytope. \par
On the other hand, if $G$ is a compact group of isometries of a Riemannian
manifold and if the $G$-action is polar, then we can associate to any
section $\Si$ a finite group $W$, also called Weyl group, which is defined
as $W := N_G(\Si)/Z_G(\Si)$; here $N_G(\Si)$ (resp.~$Z_G(\Si)$) denotes
the subgroup of $G$ given by those elements which map $\Si$ onto itself
(resp.~act trivially on $\Si$). It is well known that the intersection of
a $G$-orbit with $\Si$ coincides with a $W$-orbit and we can identify
the two orbit spaces $M/G$ and $\Si/W$; for general properties of the
Weyl group we refer to [PT1], [PTh], [Co].\par
>From Theorem 1.1 and Theorem 1.4
we obtain the following
\bigskip
\proclaim{Corollary 1.5} Let $M$ be an irreducible homogeneous K\"ahler
manifold which is acted on polarly by a compact Lie group of isometries
$G$. Then
\roster
\item the cohomogeneity of the $G$-action does not exceed the rank of $G$;
moreover the complexified group $G^\C$ and any Borel subgroup of $G^\C$
act on $M$ with finitely many orbits;
\item if $\Phi: M \to \g^*$ denotes the moment map, $\Si$ is any section
for the $G$-action and $W$ its Weyl group, then the moment map induces a
homeomorphism between $\Si/W$ and the convex polytope
$\Phi(M)\cap \t^*_+$.\endroster
\endproclaim

To prove the corollary one only has to observe that the action of $G$ is
Poisson since the full group of isometries of $M$ acts in a
Poisson fashion: indeed any flag manifold is a coadjoint orbit whose
moment map is simply the inclusion. \par
We end this section pointing out one common property of polar and
coisotropic linear actions of a compact Lie group $G$ which will be
used in section four. It is known that a polar linear
action of $G$ on a (real) vector space does not contain two
distinct irreducible equivalent $G$-submodules (see [HPTT], Proposition
2.10); if $G$ is realized
as a subgroup of $\U(N)$ for some $N \geq 1$ and if the corresponding action
on $\C^N$ is coisotropic with respect to the standard symplectic form of
$\C^N$, then every irreducible complex $G$-module appears in $\C^N$ with
multiplicity at most one, i.e. every linear coisotropic action of $G$ is
{\it
multiplicity free\/} (see for instance [HW], Theorem 1, p.~275). We also
note
here that the restriction of a linear coisotropic action to a complex
$G$-submodule is still coisotropic (see the Equivariant Mapping Lemma in
[HW], p.~272); as a consequence, a linear coisotropic action has no
nontrivial 
fixed vectors, in contrast with the
polar actions.\par
The following result states that a coisotropic
linear action is multiplicity free also from the \lq real\rq\ point of view.
\medskip
\noindent{\bf Proposition 1.6.} {\it Let $G$ be a compact subgroup of}
$\U(N)$ {\it for
some $N \geq 1$. If the corresponding action of $G$ on $\C^N$ is coisotropic
with respect to the standard symplectic form on $\C^N$, then the real
$G$-module $\C^N \simeq \R^{2N}$ does not contain two distinct equivalent
$G$-submodules. }
\medskip
The proof of this proposition will be given in the next section.

\bigskip

\noindent{\bf \S 2. Proofs of Theorem 1.1 and Proposition 1.6}
\bigskip
Let $M$ be a compact, irreducible homogeneous K\"ahler
manifold with K\"ahler metric $g$ and complex structure $J$; if we denote
by $L$ the
identity component of the full isometry group of ($M,g$), then $L$ preserves
the complex structure $J$ (see the introduction) and
we can represent $M$ as $M = L/K$, where $K$ is the centralizer
in $L$ of a torus. Note again that the irreducibility of ($M,g$) implies
that
$L$ is simple.\par
We consider a compact connected Lie subgroup $G\subset L$ of
holomorphic isometries
and we shall suppose that the $G$-action on $M$ is polar
with a section $\Si$. \par
We saw in the introduction that a submanifold of a K\"ahler manifold
is coisotropic if and only if all of its normal spaces are totally real.
The proof of Theorem 1.1 is therefore a direct consequence of the
following proposition.
\medskip
\proclaim {Proposition 2.1} The  section $\Sigma$ is totally
real.\endproclaim
\medskip
The proof of this proposition  requires some lemmas. \par
We assume that the cohomogeneity of the action is at least two since the
claim in the proposition is trivially true for
cohomogeneity one actions. First of all, if $\Si$
is a section, we may define the complex subspace
$$\H_p = T_p\Si\cap JT_p\Si\tag 2.1$$
for $p\in\Sigma$.
We next prove that $\H$ is a complex subbundle over $\Si$. It suffices to
show that if $p$ and $q$ are two
points in $\Si$ and $\sigma$ a path in $\Si$ joining them, then parallel
transport along $\sigma$ sends $\H_p$ to $\H_q$. This is
clear since $\Si$ is totally geodesic and parallel transport in $M$ commutes
with $J$. An important consequence of this argument that we will use later
is that $\H$ is a parallel subbundle of $T(\Si)$.\par

We can now define a complex subbundle $\H$ of
$T(M)|_{M_{\text{reg}}}$ by setting
 $\H_q = T_q\Si(q)\cap JT_q\Si(q)$ for $q\in M_{\text{reg}}$, where $\Si(q)$
denotes the
unique section through $q$ and $M_{\text{reg}}$ denotes the set of regular
points of the
$G$-action.
\par
Our aim is to prove that $\H_p=\{0\}$ for all $p$. In order
to do this, we will first show that $\H$ extends to a smooth integrable
subbundle of the
whole $T(M)$, whose orthogonal complement is also integrable, and then we
will
use a result on foliations to prove our claim.\par
\medskip
\proclaim{Lemma 2.2} The subbundle $\H$ can be extended to a differentiable
subbundle of the tangent bundle $T(M)$.\endproclaim
\medskip
\demo{Proof} We consider a nonregular point $q\in M$ and fix
two sections $\Si,\Si'$ passing through $q$. If we denote by
$\H_q(\Si)$ and $\H_q(\Si')$ the complex subspaces of $T_q\Si$ and
$T_q\Si'$ respectively, defined as in (2.1), we first must prove that they
coincide and then that the extension of the bundle to the whole
manifold $M$ is differentiable.

We will use the fact that the isotropy representation of $G_q$ is
polar (see [PT1,2]).
We start proving the following
\medskip
\proclaim{Sublemma 2.3} The subspace $\H_q(\Si)$ is contained in
the intersection of all singular hyperplanes of $T_q\Si$ with respect to the
isotropy representation of $G_p$.\endproclaim
\medskip
\demo{Proof} Indeed, if $Y\subset T_q\Si$ is a singular hyperplane,
then there exists an element $g\in G_q$ which leaves $\Si$ invariant and
such that $dg_q$ is a reflection in $Y$. Since $g$ is holomorphic,
it leaves the maximal holomorphic subspace of $T_q\Si$ invariant; moreover,
since the fixed point set $Y$ has real codimension one, it must contain
$\H_q(\Si)$.\qed\enddemo
\medskip
\proclaim {Sublemma 2.4} Let $G$ be a compact Lie group acting linearly
and polarly on a vector space $V$. If $\Si$ and $\Si'$ are two different
sections
with nontrivial intersection, then the intersection $\Si\cap\Si'$ contains
the intersection of all singular hyperplanes of $\Si$ (and of $\Si'$).
\endproclaim
\medskip
\demo{Proof} We use induction on the dimension of $V$, the case when
$\dim V=2$ being trivially true. First of all, we may suppose that the
fixed point set $V^G$ is equal to $\{0\}$. Indeed, if $V^G\not=\{0\}$,
we decompose $V=V^G\oplus W$, with $W=(V^G)^\perp$ and we
note that any section splits as the sum of $V^G$ plus a section for the
$G$-action in $W$; using the induction hypothesis
on $W$, we get our claim. So, we suppose that $V^G=\{0\}$. If we fix a
nonzero vector
$v\in\Si\cap\Si'$, we clearly have that $v$ is singular for the
$G$-action
and moreover the orbit $Gv$ has positive dimension. Therefore the isotropy
$G_v$ acts polarly on the normal space $N_v$ to the orbit $Gv$ and, since
$\dim N_v<\dim V$, by the induction hypothesis, we have that $\Si\cap\Si'$
contains the
intersection of all singular hyperplanes in $\Si$ passing through $v$.
A fortiori, it
then follows that $\Si\cap\Si'$ contains the intersection of all singular
hyperplanes in $\Si$.\qed
\enddemo
\medskip
We now extend the bundle $\H$ to the singular points. We fix a singular
point
$q$ and consider two sections $\Si$, $\Si'$ through $q$. They correspond to
two linear sections, also denoted by $\Si$, $\Si'$, in the normal spaces
$N_q$ of the orbit $Gq$ for the polar action of the isotropy subgroup $G_q$.
The $K$-cycles of Bott and Samelson ([BS]) can be used to prove
 that there is a sequence of linear sections $\Si_0=\Si,\dots,\Si_l=\Si'$
such that $\Si_i\cap\Si_{i+1}$ is nontrivial since we are assuming
cohomogeneity greater than one (see the proof of Lemma 1B.3 in [PTh]). From
Sublemma 2.4, we
see that $\Si_i\cap\Si_{i+1}$ contains the intersection of all singular
hyperplanes of $\Si_i$ and therefore, by Sublemma 2.3, it contains
$\H_q(\Si_i)$; similarly, it contains $\H_q(\Si_{i+1})$. Since
$\H_q(\Si_i)$ is the maximal complex subspace of $\Si_i$, it follows that
$\H_q(\Si_{i+1})\subseteq \H_q(\Si_i)$ and we
also have the opposite
inclusion,
so that $\H_q(\Si_i)=\H_q(\Si_{i+1})$. Therefore $\H_q(\Si)=\H_q(\Si')$.\par
We are left with proving that the extension is
differentiable.\par
We fix a singular point $q \in M$ and a basis $\{v_1,\dots,v_l\}$ of $\H_q$;
we will show that there exist smooth vector fields $v_1^*,\dots,v_l^*$
on some suitable neighborhood $U$ of $q$, which extend
$v_1,\dots,v_l$ respectively and span $\H_y$ for all $y \in U$. \par
We denote by $N_q$ the normal space $T_q(G\cdot q)^\perp$ and choose
a real positive number $r > 0$ so that $\exp_q|_{B_r} : B_r \to M$
is a diffeomorphism onto its image, where $B_r := \{v \in N_q;\ ||v|| <
r\}$.
We set ${\Cal N}_r := \exp_q(B_r)$ and note that ${\Cal N}_r = \bigcup
{\Cal N}_r\cap \Sigma$, where the union is taken over all sections $\Si$
passing through $q$. If $y \in {\Cal N}_r$, then $y = \exp_q(Y)$ for some
$Y \in B_r$ and for $i = 1,\dots,l$ we define ${\hat v}_i(y)$ to be the
parallel transport of $v_i$ along the geodesic $\exp_q(tY)$, $t \in [0,1]$;
note that if $\Sigma$ is any section through $q$ with $Y \in T_q\Si$, then
$y \in \Si$ and ${\hat v}_i(y) \in \H_y \subseteq T_y\Si \subset
T_y{\Cal N}_r$, since $\H$ is parallel along the sections. \par
Let now $\g_q$ be the Lie algebra of the isotropy subgroup $G_q$ and fix
a complement $\m_q$ of $\g_q$ in $\g$ so that $\g = \g_q + \m_q$; it is
clear that, if $r$ is sufficiently small, there exists a neighborhood
$V$ of $O\in \m_q$ so that the map
$$\eqalign{ V \times {\Cal N}_r & \to M\cr
(X,y) & \mapsto \exp^G(X)\cdot y}$$
is a diffeomorphism onto an open neighborhood $U$ of $q$ in $M$. If $x \in
U$,
say $x = \exp^G(X)\cdot y$ with $X \in V$ and $y \in {\Cal N}_r$, we
define $v_i^*(x) = \exp^G(X)_*({\hat v}_i(y)) \in
\H_x$ for all $i = 1,\dots,l$ and this proves our claim.\par
\bigskip
We now prove the following
\proclaim {Sublemma 2.5} The distribution given by $\H$ is integrable and
totally geodesic; the distribution given by $\H^\perp$
is integrable.\endproclaim
\demo{Proof} It is enough to prove our claims on the
regular set
$M_{\text{reg}}$ which is dense in $M$. \par
The  claims about $\H$ follow from the fact that $\H$ is a parallel
subbundle of $T(\Si)$ for all sections $\Si$. \par
We now put $k=\dim \Si-\dim \H$ and $m$ the dimension of a regular
$G$-orbit.
Then $\dim \H^\perp=k+m$.
We pick a regular point $q$ and consider the section $\Si$ through $q$;
we then
fix a set $\{Z_1,\dots,Z_k\}$ of smooth independent sections of
$\H^\perp$ which are tangent to $\Si$
over a suitable neighborhood $V\cap\Si$, where $V$ is an open neighborhood
of $q$ in $M_{\text{reg}}$. We may
extend the sections $\{Z_1,\dots,Z_k\}$ to vector fields on a neighborhood
of $q$ by using the $G$-action, since the slice representation of $G_p$
is trivial for every $p\in V$. On some suitable open neighborhood $V'$ of
$q$ in 
$V$, we can choose Killing vector fields $\{X_1,\dots,X_m\}$, induced by the
$G$-action, 
which span the tangent space to
every orbit $Gp$ for $p\in V'$; the vector fields
$\{Z_1,\dots,Z_k,X_1,\dots,X_m\}$ are therefore linearly independent
and span $\H^\perp_p$ for all $p \in V'$. We now prove that $\H^\perp$
is integrable. First of all we note that, by the definition of $Z_i$
the Lie derivative ${\Cal L}_{X_j}Z_i = 0$ and therefore $[Z_i,X_j] = 0$
for all $i = 1,\dots, k,\ j = 1,\dots,m$; moreover, $[X_i,X_j]$ is tangent
to the $G$-orbits and is therefore a
section of $\H^\perp$ for all $i,j = 1,\dots,m$. It is left to prove that
$[Z_i,Z_j]$ are sections of $\H^\perp$ for all $i,j = 1,\dots,k$.\par
Given any section $X \in \Gamma(\H)$ on $V'$, we have
$$g([Z_i,Z_j],X) = g(Z_i,\nabla_{Z_j}X) - g(Z_j,\nabla_{Z_i}X),$$
where $\nabla$ denotes the Levi-Civita connection of the metric $g$. If
$x \in V'$ and $\Si(x)$ denotes the (unique) section through $x$, then
$Z_i,Z_j$ and $X$ are tangent to $\Si(x)$ and $\nabla_{Z_i}X,\nabla_{Z_j}X$
are sections of $\H$, since $\H$ is a parallel subbundle of $T(\Si(x))$;
therefore $g([Z_i,Z_j],X) = 0$ and our claim follows.  \qed\enddemo
\enddemo
\medskip
\proclaim {Lemma 2.6} We have $\H_p = \{0\}$ for all $p\in M$.\endproclaim
\demo{Proof} In [BH]  the following result is proved: given a totally
geodesic foliation $\Cal F$ on a compact Riemannian manifold ($M,g$) such
that
the distribution given by  $T{\Cal F}^\perp$ is integrable,
then the universal covering manifold
$\tilde M$ splits topologically as the product of two manifolds that
are the universal coverings of leaves of the two distributions $\Cal F$ and
$\Cal F^\perp$. In our situation, the manifold $M$ is simply connected;
moreover, it is known (see [On], p.~279 and [Sh]) that
the cohomology ring $H^*(M,\Bbb R)$ is indecomposable, i.e. can not be
written as a tensor product $A\otimes B$, where $A,B$ are differential
graded algebras of dimension bigger than two. This implies that the leaves
of the distribution $\H$ are zero dimensional.\qed\enddemo
\medskip
\demo {Proof of Proposition 2.1} Since the manifold $M$ is compact, the
complexified Lie group $G^\C$ acts on $M$. Moreover, we can consider an
Iwasawa decomposition $\g^\C = \g+\s$, where $\s$ is a solvable Lie
subalgebra.
 We recall that the manifold $M = L/K$
can be $L^\C$-equivariantly embedded into some complex projective space.
If $S$ denotes the Lie subgroup of $G^\C$ with Lie algebra $\s$, then
by the Borel Fixed Point Theorem, $S$ has a fixed point in $M$ and therefore
there exists a point $q\in M$ such that the $G^\C$-orbit $G^\C q$ coincides
with the $G$-orbit $Gq$, which is therefore complex. We denote by $N_q$ the
normal space of the orbit $Gq$ at $q$; then $N_q$ is a complex subspace of
$T_qM$
and the isotropy $G_q$ acts on it polarly with a section $\Si$.\par
First of all, we note that the fixed point set $N_q^{G_q}$ is equal to
$\{0\}$, since otherwise it would be a complex subspace contained in a
section,
contradicting Lemma 2.6. Therefore the $G_q$-module $N_q$ splits as
a sum of irreducible, nontrivial complex submodules; these are also
irreducible over the
reals since otherwise $N_q$ would contain two equivalent real submodules,
which 
is not possible by [Da], see also [HPTT], p.~168.
Our claim will now follow from the following lemma.
\medskip
\noindent{\bf Lemma 2.7.} {\it Let $K$ be a compact Lie subgroup of\/}
$\O(V)$, 
{\it where
$V$ is a real vector space. If the action of $K$ on $V$ is irreducible,
polar
and leaves some complex structure $J$ on $V$ invariant, then the sections
are totally real subspaces with respect to $J$.}
\medskip
\demo{Proof} We know from Dadok's Theorem that the $K$-action is orbit
equivalent to the isotropy representation of a symmetric space $U/H$ and
that $K$ can be considered to be a subgroup of $H$.
Furthermore the representation of $K$ on $V$ is the restriction of the
isotropy representation of $U/H$. If
$K$ itself is the isotropy subgroup of a symmetric space $U/K$, then $U/K$
is Hermitian symmetric and the sections are totally real: indeed
if $\u = \k + \m$ is the corresponding Cartan decomposition, we
identify $V$ with $\m$ and the complex
structure $J$ can be represented as the adjoint action of some element
$z\in \k$; therefore, if $\Si$ is any section and $v,w \in \Si$, then
$<Jv,w> = <[z,v],w> = 0$, since $w$ belongs to $\Si$ which is orthogonal
to the tangent space of the orbit Ad($K$)-orbit through $v$.\par
If, on the other hand, there exists a symmetric pair ($U,H$) such that
$K$ is a proper subgroup of $H$ and $H,K$ have the same orbits in $V$, where
$\u = \h + V$ is the corresponding Cartan decomposition, then the pairs
($K,H$) have been completely classified in [EH] (see also [GT]). We note
here that our claim
is obviously true when the cohomogeneity is one. The
list of  such pairs ($K,H$) with cohomogeneity of the $K$-action at least
two is given in the following table, where we also
indicate the group $U$, such that ($U,H$) is a symmetric pair.
\bigskip
\eightpoint
\centerline{
\vbox{\offinterlineskip \halign {\strut\vrule \hfil \ $#$\ \hfil
&\vrule\hfil \
$#$\ \hfil  &\vrule\hfil\ $#$\
\hfil\vrule\cr
\noalign{\hrule }
\ \ \ \ {K}_{\phantom {\sum_1^N}} &
 \ {H}_{\phantom {\sum_1^N}} &
 \ {U}_{\phantom {\sum_1^N}}\cr
\noalign{\hrule depth 1 pt}  {\SO(2)\times \G2} &
 {\SO(2)\times \SO(7)}
 & \ {\SO(9)}\cr
\noalign{\hrule} {\SO(2)\times \Spin(7)} &
 {\SO(2)\times \SO(8)}
 & \ {\SO(10)}\cr
\noalign{\hrule} {\SO(3)\times \Spin(7)} &
 {\SO(3)\times \SO(8)}
 & \ {\SO(11)}\cr
\noalign{\hrule} {\SU(p)\times \SU(q),\ p\not=q} &
 {S(\U(p)\times \U(q))}
 & \ {\SU(p+q)}\cr
\noalign{\hrule} {\SU(n),\ n\ {\text {odd}}} &
 {\U(n)}
 & \ {\SO(2n)}\cr
\noalign{\hrule} {\SO(10)} &
 {\T^1\cdot \SO(10)}
 & \ {\E_6}\cr
\noalign{\hrule}}}
}
\tenpoint
\bigskip
Since the corresponding symmetric spaces $U/H$ are Hermitian symmetric
with the only exception of $H = \SO(3)\times \SO(8)$ and $K = \SO(3) \times
\Spin(7)$, it will be enough to prove that the representation space
$V = \Bbb R^{24}$ of the group $K$ (here the representation is the tensor
product of the standard representation of $\SO(3)$ on $\Bbb R^3$
and the spin representation of
$\Spin(7)$ on $\Bbb R^8$) does not carry any invariant complex structure.
If $V$ carries such a complex structure, then $K$ would have an irreducible
representation in $\C^{12}$, which would be a tensor product of two
irreducible representations of $\SO(3)$ and $\Spin(7)$ respectively; this
contradicts
the fact that the least dimensional irreducible complex representation of
$\Spin(7)$ is seven dimensional. \qed\enddemo
\bigskip
\remark{Remark} We point out here that the homogeneity of the K\"ahler
manifold $M$ was used only in the proof of Lemma 2.6 and at the beginning of
the proof of Proposition 2.1 where we use the Borel Fixed Point
Theorem.
Actually, looking at the proof of Lemma 2.6, we see
that Theorem 1.1 holds true for every compact, projective algebraic K\"ahler
manifold whose
universal covering space is not diffeomorphic to a product of two
K\"ahler manifolds of lower dimension.
\endremark
\bigskip
\demo{Proof of Proposition 1.6} We denote by $J$ the complex
structure and by $<,>$ the scalar product on the real $G$-module $\R^{2N}$.
We split $\C^N = \sum_{i=1}^k W_i$ as
the sum of mutually inequivalent complex irreducible $G$-submodules; we
recall that each $(W_i)_{\R}$, i.e. the real $G$-module obtained from $W_i$
by restricting the scalars, is either irreducible or splits as the sum of
two equivalent submodules which are interchanged by $J$. Suppose that the
$G$-module $\R^{2N}$ contains two equivalent irreducible $G$-submodules;
then we can suppose that there exist two equivalent $G$-submodules $E_1,E_2$
such that either $E_1,E_2$ are both complex or $E_1+E_2$ is complex and $J$
interchanges
$E_1,E_2$. Since
the restriction of the $G$-action to a complex
submodule is still coisotropic we can restrict ourselves to the complex
submodule $W:= E_1 + E_2$. We identify both $E_1,E_2$ with the irreducible
$G$-module $V$
endowed with the invariant scalar product $g$ and $W$ with the orthogonal
sum $V\oplus V$ so that
$<,>$ restricts to $g$ on the first factor and with $c^2 g$ on the second
one
for $c\in \R\setminus\{0\}$. If $E_1$ and $E_2$ are not complex
submodules, then $Hom(V,V)^G = \R$ (see e.g.~[BtD], p.~99)
and it is easy to see that $J$ is
given up to the sign by
$J(v,w) = (-cw,{1\over c}v)$ for $v,w\in V$; if both $E_1,E_2$ are complex
submodules, then $V$ has a complex structure $\tilde J$ and $J(v,w) =
({\tilde J}v, -{\tilde J}w)$, since $E_1\simeq {\overline{E_2}}$ as complex
modules. \par
We now construct two functions $f_1,f_2 \in C^\infty(W)^G$ which do
not commute with respect to the Poisson bracket, contradicting Theorem 1.4.
We define
$f_1(v,w) = g(v,v)$ and $f_2(v,w) = g(v,w)$
for $v,w\in V$. Then $X_{f_1}|_{(v,w)}$ is given by
$(2{\tilde J}v,0)$
if $V$ is complex and by $(0,{2\over c}v)$ if $JE_1 =
E_2$.
It is then easy to check
that $\{f_1,f_2\} = -df_2(X_{f_1})$ does not vanish identically, proving
our claim.\qed
\enddemo
\bigskip
\noindent{\bf \S 3. New examples of polar actions}
\bigskip
In this section we first prove Theorem 1.2 and then give new examples of
polar
actions.
We also show that there is a coisotropic action which is not polar.
\medskip
\noindent{\it Proof of Theorem 1.2.}
We will denote by $J$ and $\omega$ the complex structure and the
K\"ahler form of $M$ respectively and denote the torus $T^n$ by $G$. The
commutativity of
$G$ implies that the
action has only one principal isotropy group which is contained in all other
isotropy groups.  It follows that the principal isotropy group is trivial
since
the action is effective. Hence the principal orbits are $n$-dimensional.
Since the Euler characteristic of $M$ is
positive, the $G$-action has a fixed point and it is therefore Poisson
by a result of Frankel [Fr]; this together with the commutativity of
$G$ implies that $\omega(X,Y) = 0$ for all $X,Y \in
T_p(Gp)$ and therefore $JT_p(Gp)=N_p(Gp)$ for every regular point
$p$, where $N_p(Gp)$ denotes the normal space of $Gp$ at
$p$. We denote by $\Phi$ the corresponding moment map, which separates
$G$-orbits by Theorem 1.4.\par
Let $P=\Phi(M)$ be the moment polytope of the
$T^n$-space $M$. Delzant constructs in [De] a \lq canonical\rq\ smooth
projective variety $M_P$ that is $G$-equivariantly symplectomorphic
to $M$. The symplectic manifold $M_P$ is naturally equipped with
a K\"ahler metric $g_P$, a complex structure $J_P$ and a symplectic form
$\omega_P$; moreover there exists an antiholomorphic involution $\tau_P$
on $M_P$ such that $\tau_Pg = g^{-1}\tau_P$ for all $g \in G$ and
$\Phi\tau_P = \Phi$. We denote by $Q_P$ the fixed point set of $\tau_P$. The
compact, properly embedded submanifold $Q_P$ is not empty: indeed, if $z\in
M_P^G$, the fixed point set of $G$ in $M_P$, then $\tau_P z \in M_P^G$ and
$\Phi\tau_P z =\Phi z$. Since $\Phi$ separates orbits,
 we have $\Phi^{-1}\Phi w = w$ for all $w \in
M_P^G$. Hence $\tau_P z = z$.\par
It follows from [Du] that $\Phi(Q_P) = \Phi(M_P)$ since $Q_P$ is not empty.
Hence $Q_P$ contains
regular points; if $q \in Q_P$ is regular, then we have the orthogonal
splitting  $T_qM_P = T_q(Gq) \oplus T_qQ_P$ (with respect to $g_P$) and
$J_P(T_qQ_P) =
T_q(Gq)$.\par
Abreu shows in [Ab] that the complex manifolds ($M,J$) and
($M_P,J_P$) are $G$-equivariantly biholomorphic. If $\phi$ is such a
biholomorphism,
then $\tau := \phi^{-1}\tau_P\phi$ is an antiholomorphic involution of $M$
(which is not necessarily isometric). If $Q$ denotes the fixed point set of
$\tau$,
namely $Q := \phi^{-1}(Q_P)$, then for every regular point $q\in Q$ we have
$T_qQ = JT_q(Gq) = N_q(Gq)$. This shows that a connected component $\Si$
of $Q$ is an integral submanifold for the normal distribution to the regular
orbits. It is a well known fact that $\Si$ meets every $G$-orbit; moreover
$\Si$ intersects the regular and hence all $G$-orbits orthogonally.  \qed
\medskip

We will now give  examples of  torus actions, in which we can explicitly
describe the sections.
\medskip
\noindent{\bf Example 1.} We consider
  the complex projective space $\PC$
together with the polar action of the
$n$-dimensional torus $T^n$, which is simply the maximal torus in the group
$\SU(n+1)$; 
the regular orbits are $n$-dimensional tori and any section is a totally
geodesic real projective space ${\Bbb P}_n({\Bbb R})$. \par
Our aim is to produce new examples of polar actions on
nonhomogeneous K\"ahler manifolds, which are not hyperpolar.
In order to do this, we choose a $T^n$-fixed point $p \in \PC$ and we denote
by $M_1$ 
the complex manifold which is obtained blowing up $\PC$ in $p$; note that
$M_1$ is a projective variety and that $M_1$ is not homogeneous
with respect to the full group of holomorphic automorphisms (hence with
respect to the full group of isometries of any K\"ahler metric compatible
with 
the given complex structure), see [Bl]. The $T^n$-action lifts to
a  holomorphic action on $M_1$; moreover, since $M_1$ is K\"ahler and $T^n$
is compact, we can 
find a $T^n$-invariant K\"ahler metric $g$ on $M_1$ by averaging; we remark
that the metric $g$ has
{\it a priori\/} no relation with the Fubini-Study metric on $\PC$. It
follows 
from Theorem 1.2 that the action of $T^n$ on $M_1$ is polar. In the
following
we will determine the section of the action.\par
We
denote by $\pi: M_1 \to \PC$
the holomorphic projection.
We fix a section $\Si$ in $\PC$ and denote by $\Si_1$ the closed
submanifold of $M_1$ which
coincides with $\Si$ outside $\{p\}$ and has the property that $\Si_1 \cap
\pi^{-1}(p)$ is ${\Bbb P}_\C(T_p\Si)$,
where we have identified $\pi^{-1}(p)$ with ${\Bbb P}_\C(T_p\PC)$. It is
clear that $\Si_1$ meets every
$T^n$-orbit. We are now going to show that $\Si_1$ is a section for the
$T^n$-action with respect to the
K\"ahler metric
 $g$ and we will prove this in the dense open set of
regular points. \par
First of all we remark that the $T^n$-action on $M_1$ is a Poisson action,
since the manifold $M_1$ is simply connected (see e.g. [GS3]);
moreover the $T^n$-action is
coisotropic by Theorem 1.4, (v). We
denote by $J$ the complex structure on $M_1$ and
fix a regular point $y \in M_1\setminus \pi^{-1}(p)$.  If we denote by $N_y$
the normal space to
the  orbit $T^n\cdot y$, then $JN_y = T_y(T^n\cdot y)$ since the action is
coisotropic. On the other hand,
$J(T_y\Si_1)$ also coincides with $T_y(T^n\cdot y)$,
so that $N_y = T_y(\Si_1)$ and therefore $\Si_1$ is a section.
We remark that the section $\Si_1$, as a differentiable manifold, is the
connected sum $\PR \# \PR$. \par
Since the $T^n$-action on $M_1$ has fixed points in $\pi^{-1}(p)$,
we can iterate this construction blowing up $M_1$ in one of
these fixed points.
In this way, we get a sequence of compact
K\"ahler manifolds $M_k$ on which the group $T^n$ acts polarly with the
sections being the connected sum of
$k + 1$ copies of $\PR$. If $n \geq 3$, then $\pi_1(\Si_1)$ is a free group
$\Bbb Z_2 * \Bbb Z_2 * \ldots
* \Bbb Z_2$ ($k + 1$ times). Therefore it can not carry a flat metric by the
Bieberbach Theorem (see [Wo]).
\medskip
\remark{Remark} We point out that the lifted action of $T^n$ on $M_k$ is
polar with respect to {\it any\/} invariant
K\"ahler metric. We also remark that this construction yields polar actions
on four dimensional K\"ahler manifolds with
any nonorientable compact surface as a section.\endremark
\medskip
We would now like to give examples of polar actions of nonabelian Lie
groups using the method of blow-ups. \medskip
\noindent {\bf Example 2.} We first describe the general principle which
will allow us to obtain new examples. We consider
a compact Lie group $G$ acting polarly on a homogeneous K\"ahler manifold
$M$ with at least one fixed point $p$; we denote by $K$
the  isotropy subgroup of a regular point $q\in M$. We recall the
following well-known facts on moment maps (see [GS3] and also
[HW]).  Denote by $\Phi: M \to \g^*$ the moment map with respect to the
given Poisson-K\"ahler action of $G$ and
recall that
$\Phi$ maps $G$-orbits onto coadjoint orbits in $\g^*$; moreover, for any
regular point $y\in M$, the isotropy subgroup
$H_y:= G_{\Phi(y)}$ contains $G_y$ and, at the level of Lie algebras, we may
decompose $\h_y$ (with respect to a fixed biinvariant
scalar product on $\g$) as
$$\h_y = \g_y + \a_y,$$
where $\a_y = \g_y^\perp \cap \h_y$ is abelian and centralizes $\g_y$.  If
we denote by
$N_y$ the normal space to the regular orbit $G\cdot y$ and by $J$ the
complex structure on $M$, then
$$J(N_y) = T_y(H_y\cdot y) \subset T_y(G\cdot y),$$
since the $G$-action on $M$ is coisotropic by Theorem 1.1. It follows that
the
dimension of $\a_y$ is equal to the cohomogeneity of the action of
$G$.
The tangent space
to the orbit $G\cdot y$ can be identified
with the complement $\g_y^\perp$ and we have the decomposition
$$\g = \g_y + \a_y + \m_y,$$
where $\m_y$ corresponds to the maximal complex subspace of the tangent
space $T_y(G\cdot y)$. In the spirit of S. Lie and
F. Engel ([LE], p.~501 ff), we call the
$G$-action
$\C$-{\it asystatic\/} if for a regular orbit $G/K$ there is no
nontrivial complex subspace of the tangent space $T_{[eK]}(G/K)$ which is
left fixed by the isotropy representation of $K$. (For a
modern treatment of asystatic actions, see [AA] and [PT1]). It follows
 that if a polar actions as above is also
$\C$-asystatic, then the space of
all tangent vectors to a regular orbit which are fixed under the isotropy
representation has dimension equal to the
cohomogeneity of the $G$-action.
\par
We now assume that the $G$-action is $\C$-asystatic and we consider the
blow-up $M_1$ of $M$ in the fixed point $q$.  If $\Si$
denotes a section for the $G$-action on $M$, then the submanifold $\Si_1$ of
$M_1$ constructed exactly as in Example 1 is a
closed submanifold which meets every $G$-orbit in $M_1$. If we now endow
$M_1$ with {\it any\/} $G$-invariant K\"ahler metric $g$,
we show that $\Si_1$ is a section, namely that it meets every $G$-orbit
orthogonally. Indeed, the lifted action of $G$
on $M_1$ is Poisson, since $M_1$ is simply connected; moreover the
$G$-action is coisotropic because $M_1$ is a projective
algebraic variety, $G^\C$ acts on $M$, hence on $M_1$, with an open orbit
and the same is true for a Borel subgroup of $G^\C$
(see Theorem 1.4). We fix a regular point $y \in M_1$ and
we denote by $\tilde N_y$
the normal space to $T_y(G\cdot y)$ with respect to
$g$, so that $J{\tilde N_y} \subseteq T_y(G\cdot y)$; since the action of
$G$ on $M$ is also coisotropic, we have that
$J(T_y\Si_1) \subseteq T_y(G\cdot y)$ and we claim that $J(\tilde N_y) =
J(T_y\Si_1)$. Both subspaces $J(\tilde N_y)$ and
$J(T_y\Si_1)$ are fixed under the isotropy representation of $G_y$ and have
the same dimension, namely the cohomogeneity of $G$
and therefore, by
the assumption that the action is $\C$-asystatic, they
coincide.\par
As an application, we consider the quadric $M := Q_n =
\SO(n+2)/\SO(2) \times \SO(n)$. The isotropy subgroup
$G = \SO(2) \times \SO(n)$ acts polarly on $Q_n$ with sections being two
dimensional flat tori. It is easy to check that a regular
isotropy subgroup $K$ is given by $K = {\Bbb Z_2}\times \SO(n-2)$
 and therefore the isotropy
representation splits as $2\Bbb R + 2\Bbb R^{n-2}$,
where $\Bbb R$ is a trivial submodule, while $K$ acts on $\Bbb R^{n-2}$ in
the standard way; therefore the $G$-action is
$\C$-asystatic. The blown up manifold $M_1$ is acted on polarly by $G$ with
sections being the connected sum $T^2 \# {\Bbb
P}_2({\Bbb R})$. Note that $M_1$ is not homogeneous, as already pointed out
in Example 1.\medskip

\noindent{\bf Example 3.} We will now give  examples of
coisotropic action on complex projective spaces
that are not polar. Let $G$ denote the compact Lie group $\U(2)\times
\Sp(n)$, 
$n\ge 2$,
acting on $\C^{4n}=\C^2\otimes \C^{2n}$ as a tensor product
of the standard
representations of $\U(2)$ and $\Sp(n)$ respectively. This representation
is not polar. (It has
cohomogeneity three and is shown to have taut orbits in [GT], but these
properties will not
be important here.)
We consider the complexification $G^\C=\text{GL}(2,\C)\times
\text{Sp}(n,\C)$ of
$G$ and the corresponding representation on $\C^{4n}$ which is multiplicity
free by
the classification of such irreducible representations of reductive
noncompact complex
Lie groups in [Ka]. This means that $G^\C$ acts on $\C^{4n}$ with an open
orbit as well 
as every Borel subgroup of $G^\C$. The induced action of $G^\C$ on
${\Bbb{P}_{4n-1}(\Bbb C)}$
clearly also has the property that every Borel subgroup of $G^\C$ acts on
${\Bbb{P}_{4n-1}(\Bbb C)}$
with an open orbit. This implies by Theorem 1.4 that the
action of $G$ on
${\Bbb{P}_{4n-1}(\Bbb C)}$
is coisotropic. This action is not polar since otherwise the action of
$\T^1\times G$ on $\C^{4n}$,
where $\T^1$ acts as multiplication by scalars, would have to be polar by
[PTh]; since the center of $G$ already acts on $\C^{4n}$ as
multiplication by scalars by Schur Lemma, $\T^1\times G$ and $G$ have the
same orbits and therefore the $G$-action on
${\Bbb{P}_{4n-1}(\Bbb C)}$ is not polar.

It is not difficult to classify coisotropic
actions on complex
projective spaces.
 If $G$ is a compact subgroup
of $\U(n+1)$, then the action of $G$ on $\PC$ is coisotropic if and only if
the Borel subgroups of
 ${\Bbb C}^*\times G^{\Bbb C}$ have dense orbits in ${\Bbb C}^{n+1}$ which
in
turn is equivalent to the action of ${\Bbb C}^*\times G^{\Bbb C}$ being
multiplicity free on
${\Bbb C}^{n+1}$ (see e.g. [Kr], p.~199).
One can now use the classification of multiplicity free representations
of reductive groups that can be found in [Ka] for the irreducible case and
in [BR]
and [Le] for the reducible case.

\bigskip
\bigskip
\noindent{\bf \S 4. Polar and coisotropic actions on quadrics}
\bigskip
The aim of this section is to prove Theorem 1.3, or more precisely to
classify
coisotropic and polar actions on the Hermitian
symmetric space $Q_{n-2} = \SO(n)/\SO(2)\times \SO(n-2)$ ($n\ge 5$), which
can be
represented
as the $(n-2)$-dimensional complex quadric in ${\Bbb P}_{n-1}(\C)$, i.e.,
the
set of all totally isotropic lines in $\C^n$ with respect to the
invariant symmetric bilinear form $<,>$ given by $<z,w> = \sum_{i=1}^n
z_iw_i$, or equivalently as the Grassmannian of oriented two-planes in
$\R^{n}$.
\par
Hyperpolar actions of compact Lie groups on irreducible
compact symmetric
spaces were classified by Kollross ([Ko]). We state his result in the
special
case of complex quadrics in the following theorem.
\medskip
\noindent{\bf Theorem 4.1 (Kollross).} {\it Let} $G\subset \SO(n)$ {\it be a
compact
connected Lie subgroup acting on the symmetric space} $Q_{n-2} =
\SO(n)/\SO(2)\times \SO(n-2)$ {\it in a hyperpolar fashion. Then there
exists
a compact connected Lie subgroup $G'$ of\/} $\SO(n)$ {\it containing $G$
such
that the $G'$-action on $Q_{n-2}$ has the same orbits as the $G$-action
and $G'$ is conjugate to one of the following subgroups}
\roster
\item $G' = \SO(p)\times \SO(q) \subset \SO(p+q)$ with $p,q \geq 1$;
\item $G' = \U(n) \subset \SO(2n)$;
\item $G' = \Sp(1)\cdot\Sp(n) \subset \SO(4n)$;
\item $G' = \Spin(9) \subset \SO(16)$.
\endroster
\medskip
\remark{Remark} The actions listed in Theorem 4.1 all have
cohomogeneity one  except case (1) when
$p,q > 1$ which has cohomogeneity two. Recall that a compact
Lie subgroup of $\SO(n)$ acts transitively on $Q_{n-2}$ if and only if $G =
\SO(n)$ with the two exceptions $G = \G2$ ($n = 7$) and
$G = \Spin(7)$ ($n = 8$) as can be seen from [On], Theorem 2, p.~226.
\endremark
\bigskip
We now consider a compact connected Lie subgroup $G$ of $\SO(n)$ acting
coisotropically on $Q_{n-2}$. The
complexification $G^\C$ acts on $Q_{n-2}$ with an open orbit $\Omega$ and
there exists a Borel subgroup of $G^\C$ which has an open orbit in
$\Omega$. \par
In order to prove Theorem 1.3, we will make an extensive use of the results
by Kimelfeld ([Ki]), who classified all reductive subgroups $G^\C$ of
$\SO(n,\C)$ having an open orbit on the quadric $Q_{n-2}$. \par
We will distinguish between the
two cases when $G \subset \SO(n)$ acts irreducibly on $\R^n$ and when $G$
acts reducibly on $\R^n$.
\medskip
\noindent{\it First case: the $G$-action on $\R^n$ is irreducible\/}
\medskip
In this case the complexified action of $\g^\C$ on $\C^n$ can be
irreducible or reducible (as the sum of two submodules $W \oplus W^*$,
where $W \simeq \C^{n/2}$). These two subcases appear as separate cases
in Kimelfeld's paper (see Tables I and II in [Ki], p.~535-536) and we too
will
reproduce them separately in Tables I and II below. \par
As for notation, we will use the standard symbols $A_s,B_s,C_s,D_s$
($s \geq 1$) and $G_2, F_4$ to denote complex simple Lie algebras; by
$R(\phi)$
we will denote the irreducible representation of a simple Lie algebra with
highest weight $\phi$ and $\phi_i$ will indicate the highest weight of
the $i$th fundamental representation of a simple complex Lie algebra
(the standard representation of $\so(n)$ on $\C^n$  will be call simplest);
by
$\t_n$ we will indicate the Lie algebra of the $n$-dimensional complex
torus and $\epsilon$ will denote the faithful one dimensional representation
of $\t_1$.\par
In the following Table I we list all complex reductive Lie algebras
$\g^\C$ with the corresponding complex irreducible representation such that
the induced action on the quadric $Q$ has an open orbit $\Omega$;
moreover we will also indicate a generic stability subalgebra $\h$,
with $p$-dimensional nilradical $U_p$, so that
$\Omega = G^\C/H$, where $H$ is an appropriate subgroup of $G^\C$ having
$\h$ as Lie algebra. \par
\bigskip\eightpoint
\centerline{\vbox{\offinterlineskip \halign {\strut\vrule \hfil \ $#$\
\hfil &\vrule\hfil \
$#$\ \hfil  &\vrule \hfil\ $#$\ \hfil
&\vrule\hfil\ $#$\
\hfil\vrule\cr
\noalign{\hrule }
\ \ {n.}_{\phantom {1}} &
 \ {\g^\C}_{\phantom {\sum_1^N}} &
 \ {\text{Represen.}}
 & \ {\h}_{\phantom {\sum_1^N}}^{\phantom  {\sum_1^N}} \cr
\noalign{\hrule depth 1 pt}
{1} &
 {\so(s),s\geq 3}
 & \ {\text{Simplest}} &\ {(\so(s-1)\oplus \t_1)+ U_{s-2}} \cr
\noalign{\hrule} {2} &
 {C_s \oplus A_1, s\geq 2}
 & \ {R(\phi_1)\times R(\phi_1)} &\ {(C_{s-2} \oplus A_1\oplus \t_1)+
U_{4s-5}}
\cr
\noalign{\hrule} {3} &
 {G_2}
 & \ {R(\phi_1)} &\ {A_2} \cr
\noalign{\hrule} {4} &
 {B_3}
 & \ {R(\phi_3)} &\ {G_2}\cr
\noalign{\hrule} {5} &
 {B_4}
 & \ {R(\phi_4)} &\ {B_3} \cr
\noalign{\hrule} {6} &
 {A_2}
 & \ {R(\phi_1+\phi_2)} &\ {\t_2} \cr
\noalign{\hrule} {7} &
 {C_2}
 & \ {R(2\phi_1)} &\ {\t_2}\cr
\noalign{\hrule} {8} &
 {G_2}
 & \ {R(\phi_2)} &\ {\t_2} \cr
\noalign{\hrule} {9} &
 {A_1}
 & \ {R(4\phi_1)} &\ {0} \cr
\noalign{\hrule} {10} &
 {C_3}
 & \ {R(\phi_2)} &\ {A_1 \oplus A_1 \oplus A_1} \cr
\noalign{\hrule} {11} &
 {F_4}
 & \ {R(\phi_1)} &\ {D_4} \cr
\noalign{\hrule} {12} &
 {C_s \oplus C_2, s\geq 2}
 & \ {R(\phi_1)\times R(\phi_1)} &\ {C_{s-2} \oplus A_1 \oplus A_1} \cr
\noalign{\hrule} {13} &
 {A_1 \oplus A_1}
 & \ {R(\phi_1)\times R(3\phi_1)} &\ {0} \cr
\noalign{\hrule}}}
}\medskip
\centerline{Table I}
\tenpoint
\bigskip
We start noting that the cases 1,3 and 4 correspond to transitive actions
of the compact real form $G$, while the cases 2 and 5 are known as
cohomogeneity one actions (see Kollross' theorem). \par
In the cases 6 to 8, the group $G$ acts on
$\g$ via adjoint representation; the cohomogeneity of the $G$-action
on the corresponding quadric $Q$ is easily seen to be at least four,
hence bigger than the rank of $G$. This means that the $G$-action on $Q$
is not coisotropic by Theorem 1.4.\par
We are now going to
exclude all other cases from the above table; in order to do this, we use
the fact that a Borel subgroup of $G^\C$ has an open orbit in $G^\C/H$,
see Theorem 1.4. This immediately excludes cases 9 and 13.\par
 The following result
(see [Ki], p.~572 ff) is well known:\par
\medskip
\proclaim{Lemma 4.2} Given a semisimple complex algebraic group $U$ and a
Borel
subgroup $B$, a reductive subgroup $H$ of $U$
has an open orbit in $U/B$
if and only if the space of $H$-fixed
vectors in every irreducible
$U$-module is at most one-dimensional.\endproclaim
\medskip
Note that, if a Borel subgroup $B$ of $G^\C$ has an open orbit in $G^\C/H$,
then
$H$ has an open orbit in $G^\C/B$. \par
In the remaining cases 10 to 12, we note that these are complexifications of
isotropy
representations of rank two symmetric spaces and that $\h$ is exactly the
complexification of a regular stability subalgebra for these
representations. Since $H$ stabilizes a two dimensional subspace of
the corresponding irreducible representation of $G^\C$, the $G$-action on
the quadric is
not coisotropic by Theorem 1.4 (v) and Lemma 4.2.\par
We now come to the subcase that the complexified action of $\g^\C$ on $\C^n$
is reducible (as the sum of two submodules $W \oplus W^*$,
where $W \simeq \C^{n/2}$).
In the following Table II, we reproduce all such cases having an open orbit
on the
corresponding quadric, only indicating the
irreducible representation $W$ and the generic stability subalgebra $\h$ for
the $G^\C$-action on the corresponding quadric; we
will indicate 
by $U_p$ the $p$-dimensional
nilradical of the generic stability subalgebra $\h$.
\bigskip\eightpoint
\centerline{\vbox{\offinterlineskip \halign {\strut\vrule \hfil \ $#$\
\hfil &\vrule\hfil \
$#$\ \hfil  &\vrule \hfil\ $#$\ \hfil
&\vrule\hfil\ $#$\
\hfil\vrule\cr
\noalign{\hrule }
\ \ {n.}_{\phantom {1}} &
 \ {\g^\C}_{\phantom {\sum_1^N}} &
 \ {\text{Represen.}}
 & \ {\h}_{\phantom {\sum_1^N}}^{\phantom  {\sum_1^N}} \cr
\noalign{\hrule depth 1 pt}
{1} &
 {A_s,s\geq 2}
 & \ {R(\phi_1)} &\ {(A_{s-2} \oplus \t_1) + U_{2s-1}} \cr
\noalign{\hrule} {2} &
 {A_s \oplus \t_1, s\geq 2}
 & \ {R(\phi_1)\times \epsilon} &\ {(A_{s-2} \oplus \t_2) + U_{2s-1}}
\cr
\noalign{\hrule} {3} &
 {A_s \oplus A_1, s\geq 2}
 & \ {R(\phi_1)\times R(\phi_1)} &\ {A_{s-2} \oplus \t_1} \cr
\noalign{\hrule} {4} &
 {A_s \oplus A_1 \oplus \t_1, s\geq 2}
 & \ {R(\phi_1)\times R(\phi_1)\times \epsilon} &\ {A_{s-2} \oplus \t_2}\cr
\noalign{\hrule} {5} &
 {\so(s) \oplus \t_1}
 & \ {{\text{Simplest}}\times \epsilon} &\ {\so(s-2)} \cr
\noalign{\hrule} {6} &
 {A_1 \oplus \t_1}
 & \ {R(2\phi_1)\times \epsilon} &\ {0} \cr
\noalign{\hrule} {7} &
 {B_3 \oplus \t_1}
 & \ {R(\phi_3)\times \epsilon} &\ {A_2}\cr
\noalign{\hrule} {8} &
 {G_2 \oplus \t_1}
 & \ {R(\phi_1)\times \epsilon} &\ {A_1} \cr
\noalign{\hrule} {9} &
 {C_s, s\geq 3}
 & \ {R(\phi_1)} &\ {(C_{s-2} \oplus \t_1) + U_{4s-5}} \cr
\noalign{\hrule} {10} &
 {C_s \oplus \t_1, s\geq 3}
 & \ {R(\phi_1)\times \epsilon} &\ {(C_{s-2} \oplus \t_2) + U_{4s-5}} \cr
\noalign{\hrule} {11} &
 {A_4}
 & \ {R(\phi_2)} &\ {A_1 \oplus A_1} \cr
\noalign{\hrule} {12} &
 {A_4 \oplus \t_1}
 & \ {R(\phi_2)\times \epsilon} &\ {A_1 \oplus A_1 \oplus \t_1} \cr
\noalign{\hrule} {13} &
 {D_5}
 & \ {R(\phi_5)} &\ {A_3} \cr
\noalign{\hrule} {14} &
 {D_5 \oplus \t_1}
 & \ {R(\phi_5)\times \epsilon} &\ {A_3 \oplus \t_1} \cr
\noalign{\hrule}}}
}\medskip
\centerline{Table II}\tenpoint
\bigskip
We start noting that case 2 is of cohomogeneity one (see Theorem 4.1) and
the
same holds for case 1: indeed, the action of $\U(s)$ on
$\SO(2s)/\SO(2)\times
\SO(2s-2)$ has a complex singular orbit which is ${\Bbb P}_{s-1}(\C) =
\U(s)/
\U(1)\times \U(s-1)$; this orbit also is an orbit
of the semisimple part $\SU(s)$ and the stabilizer $\text{S}(\U(1)\times
\U(s-1))$ still acts
with cohomogeneity one on the normal space.\par
Cases 6 to 8 and 11 to 14 can be ruled out simply because the dimension
of a Borel subgroup of $G^\C$ is strictly less than the dimension of
the open orbit $G^\C/H$. \par
In cases 3 and 4 we consider the action of the compact group $G$ acting
on the quadric $Q:= Q_{4s+2}$ and we will prove that in case 4 this
$G$-action is
not coisotropic; by Theorem 1.4 this implies that the $G$-action in
case 3 is also not coisotropic. We may take $G$ as $S(\U(2)\times \U(s+1))$
acting on $\C^2\otimes\C^{s+1}\simeq \R^{4s+4}$; we then fix two unit
vectors
$e_1\in \C^2$, $v \in \C^{s+1}$ and we take the oriented two
plane $\pi$ in $\R^{4s+4}$ given by the complex line $\C\cdot e_1\otimes v$.
It is easy to see that $G_{\pi}$ is isomorphic to $S(\U(1)\times \U(1)
\times \U(1) \times \U(s))$, so that the orbit $G/G_{\pi}$ is complex and
biholomorphic to $S^2\times {\Bbb P}_s(\C)$. Now, if $V$ denotes the
subspace
$(\C\cdot v)^\perp\subset \C^{s+1}$ and $e_2$ is a unit vector
generating $(\C\cdot e_1)^\perp\subset \C^2$, then we have
$$T_\pi Q = Hom(\pi,\pi^\perp) = \pi^*\otimes_{\R} [(e_1\otimes V)_{\R} +
(e_2\otimes V)_{\R} + (\C\cdot e_2\otimes v)_{\R}],$$
which contains two pairs of equivalent
$2s$-dimensional $G_\pi$-submodules. Since the
isotropy representation of $G_\pi$ on $T_{eG_\pi}G/G_\pi$
contains only one $2s$-dimensional submodule,
 it then follows that the slice
representation of $G_\pi$ at $\pi$ contains at least two
distinct equivalent submodules.
If the $G$-action on $Q$ is coisotropic, then the slice representation
of $G_\pi$ is also coisotropic (see [HW], p.~274), contradicting Proposition
1.6.\par
Case 5 can be ruled out because the pair ($\so(s)+\t_1,\so(s-2)$) is not a
spherical pair; indeed, if it were spherical, the pair ($\so(s), \so(s-2)$)
would be spherical too, contradicting Lemma 4.2.\par
In case 9, we show that the action of $G = \Sp(s)$ is not
coisotropic.
We consider the $G$-action on ${\Bbb H}^s \simeq \C^{2s} \simeq
\R^{4s}$ and the corresponding action on the quadric $Q:= Q_{4s-2}$. We fix
an
oriented two plane $\pi$ in $\R^{4s}$, which is invariant under the complex
structure of $\C^{2s}$. Then $G_\pi = T^1\times \Sp(s-1)$ and the orbit
$G\cdot \pi$ is diffeomorphic to a complex projective space
${\Bbb{P}_{2s-1}(\Bbb C)}$. The tangent space $T_\pi Q$ can be identified
with $Hom(\pi,\pi^\perp)$ and we have the orthogonal
splitting $\pi^\perp = \hat\pi
+ W$, where $\pi + \hat\pi$ is the quaternionic span of $\pi$, while
$W\simeq
\R^{4(s-1)}$ is the standard representation space of $\Sp(s-1)$.\par
If we denote by $\rho$ the representation of the center $T^1$ of $G_\pi$
on $\pi$, then $T^1$ acts on $\hat \pi$ by means of $\rho^{-1}$; the
tangent space $T_\pi Q$ splits as a $G_\pi$ module as $2\R + \R^2 + 2W$,
where $\R$ denotes the trivial module, $\R^2$ is acted on by $T^1$ by means
of
$\rho^2$, while $W$ is the standard representation space of $G_\pi$. Since
the isotropy representation of $G_\pi$ on the orbit $G\cdot\pi$ is $\R^2 +
W$,
we see that the normal space $N_\pi$ is isomorphic to $2\R + W$ as a
$G_\pi$-module; it follows that a regular isotropy subgroup $K$ of the
$G$-action on $Q$ is isomorphic to $T^1\times \Sp(s-2)$ of corank one in
$G$,
while the $G$-cohomogeneity on $Q$ is three. This proves that the $G$-action
on $Q$ is not coisotropic.\par
In case 10 the group $G = T^1\times \Sp(s)$ acts on $\C^{2s} \simeq \R^{4s}$
in a standard way; for notational reasons we will denote by $Z$ the one
dimensional center of $G$. We will show that the $G$-action on $Q:=
Q_{4s-2}$
is coisotropic but not polar. As in the previous case, we fix an oriented
two plane $\pi \subset \R^{4s}$ which is $Z$-invariant; the isotropy
subgroup $G_\pi$ is $Z\times T^1\times \Sp(s-1)$ and the orbit ${\Cal O}:=
G\cdot\pi$ is a complex projective space ${\Bbb{P}_{2s-1}(\Bbb C)}$. If
we compute the slice representation along the same lines as in case 9, we
see
that the normal space $N_\pi$ splits as a $G_\pi$-module as
$\R^2 + \R^{4(s-1)}$, where $G_\pi$ acts nontrivially on each factor (this
action can be shown to be actually polar). A regular isotropy subgroup
$K$ is isomorphic to $T^1\times \Sp(s-2)$ of corank two in $G$ and the
cohomogeneity is two, hence the $G$-action on $Q$ is coisotropic.\par
We now show that this action is not polar. It is not difficult to see that
there exists a totally geodesic submanifold $F \subset Q$ with $T_\pi F =
N_\pi$ such that $F$ is isometric to the complex
projective space ${\Bbb{P}_{2s-1}(\Bbb C)}$ after an appropriate rescaling
of the metrics.
Moreover we observe that the
group $G$ lies in $\U(2s)$, which acts on $Q$ with cohomogeneity one; the
two singular $\U(2s)$-orbits are $\Cal O$ and ${\Cal O}^*:= \U(2s)\cdot
\pi^* =
G\cdot \pi^*$, where $\pi^*$ denotes the two-plane $\pi$ with reversed
orientation. The orbit $\Cal O^*$ consists actually of all points of
maximal distance $d$ from $\Cal O$ and $d$ is the diameter of $F$.
Now suppose that the $G$-action is polar. If  $\Si$ is a section
passing through $\pi$, then $\Si \subset F$ and  we know from [PTh] that
$\Si$ is
 a real projective plane.
The
submanifold
$\Si$ contains a circle $S$ of points at distance $d$ from $\pi$, hence
$S$ is contained in $\Cal O^*$; this means that $\Si$ intersects the
$G$-orbit
in a circle, hence not orthogonally; a contradiction.
\bigskip
So far we have proved the following
\bigskip
\noindent{\bf Proposition 4.3.} {\it Let $G$ be a compact connected Lie
subgroup
of\/} $SO(n)$, {\it where $n \geq 5$. Suppose the action of $G$ on $\R^n$ is
irreducible. Then the action of $G$
on $Q:= Q_{n-2}$ is nontransitive and polar
if and only if it is
hyperpolar and\/} $G = \U(n/2),SU(n/2), G = \Sp(1)\cdot \Sp(n/4)$ or $G =
\Spin(9)
\subset \SO(16)$.\par
{\it The action of $G$ on $Q$ is nontransitive and coisotropic if and
only if the action is polar or\/} $G = T^1\cdot\Sp(n/4)$\medskip
\bigskip
\noindent{\it Second case: the $G$-action on $\R^n$ is reducible\/}
\medskip
We will now suppose that the linear $G$-action on $\R^n$ is reducible.
Here the proof will be more direct than in the first case.
We
will first deal with the particular case when $G$ has a fixed nonzero
vector $e$ in $\R^n$, so that we can write $\R^n = \R\cdot e + \R^{n-1}$
(note that we are not supposing that $\R^{n-1}$ is $G$-irreducible).
We prove the following
\medskip
\noindent{\bf Proposition 4.4.} {\it Let $G$ be a compact connected Lie
group
acting on $V:= \R^n$ with nonzero fixed point set. Then}
\roster
\item "(i)" {\it the action of $G$ on $Q:= Q_{n-2}$ is coisotropic if and
only if the fixed point set $V^G$ is one dimensional and $G$ acts
transitively on the unit sphere of $(V^G)^\perp \simeq \R^{n-1}$ with
the exception of\/} $G = \Sp(k)$ {\it acting in a standard way on $\R^{4k},
n = 4k +1$.}
\item "(ii)" {\it the action of $G$ on $Q:= Q_{n-2}$ is polar if and
only if the fixed point set $V^G$ is one dimensional and} $G = \SO(n-1),
G = \G2$ {\it ($n =8$)  or} $G = \Spin(7)$ {\it  ($n = 9$); in particular
the action of $G$ is polar if and only if it is hyperpolar.
}\endroster
\medskip
\demo{Proof} We suppose that the $G$-action on $Q_{n-2}$ is
coisotropic.
We first show that $G$ must act transitively on the unit sphere of
$\R^{n-1}$.
We fix a nonzero vector $v \in \R^{n-1}$ which is regular for the
$G$-action on $\R^{n-1}$ and consider the oriented
two plane $\pi = e\w v$; we note that $G_v = G_\pi$. We may identify the
tangent space $T_\pi Q$ with
$e \w W + v \w W$, where $W\subset \R^{n-1}$ is the set of all vectors
orthogonal to $v$ in $\R^{n-1}$; moreover $T_\pi G\cdot\pi$ is given
by
$$T_\pi G\cdot\pi = \{Xe\w v + e\w Xv;\ X\in\g\} = \{e\w Xv;\ X\in\g\}.$$
We note that $\g\cdot v \subset W$ and we split orthogonally
$W = \g\cdot v + V$,
where $\dim V + 1$ is equal to the cohomogeneity of $G$ in $\R^{n-1}$.
Then 
$$T_\pi Q = T_\pi G\cdot\pi + e \w V + v \w V + v \w \g\cdot v.$$ We recall
then that the complex structure $J$ on $T_\pi Q$ is given (up to sign) by
$$J(e \w w_1 + v \w w_2) = e \w w_2 - v \w w_1,$$
for $w_1,w_2 \in W$. Therefore the subspace $e \w V + v \w V$ is a
complex subspace of $(T_\pi G\cdot\pi)^\perp$ on which $G_\pi$ acts
trivially.
The slice representation of $G_\pi$ is coisotropic
since the action of
$G$ on $Q_{n-2}$ is coisotropic (see [HW], p.~274) and can therefore not
leave a nontrivial subspace invariant. This implies that
that $e \w V + v \w V=\{0\}$, hence $V = \{0\}$ and our first
claim follows.\par
We now use the fact that $G$ induces a coisotropic, possibly transitive
action on the subquadric $Q(\R^{n-1}) =
Q_{n-3}$. An inspection of the lists of compact groups acting
transitively on spheres (see e.g. [Bo]) and quadrics (see the remark after
Theorem 4.1)
together with Proposition 4.3 shows that $G$ can be
either \par
({\it a\/}) $G= \SO(n-1)$,\; $\G2$, or $\Spin(7)$ which all
give rise to a cohomogeneity one action on $Q$ (those of $\G2$ and
$\Spin(7)$ being
orbit equivalent to the one of $\SO(n-1)$ for an appropriate
$n$), \par
({\it b\/})
$G=\SU({{n-1}\over 2})$, $\U({{n-1}\over 2})$, $\Sp(1)\cdot\Sp({{(n-1)}\over
4})$, $\T^1\cdot\Sp({{(n-1)}\over 4})$ or $\Spin(9)$. \par
It is easy to check that all groups listed in ({\it b\/}) yield coisotropic
actions; their cohomogeneity on the quadric is given by the cohomogeneity
of $G_v$ acting on $\g \cdot v$.\par
We will now show that the $G$-action on the quadric is not polar
in all the cases enumerated in ({\it b\/}).
Let now $G$ be one of the groups listed in ({\it b\/}); we already
know that the induced action of $\T^1\cdot\Sp({{(n-1)}\over 4})$ on the
subquadric $Q_{n-3}$ is not polar if $n > 5$, while, for $n = 5$, $G$
coincides with $\U(2)$, so that this case can be ruled out. Using
the same notation as above, it is easy to see that $G_v$ acts on $\g \cdot
v$ polarly and with cohomogeneity two. We choose  two  orthogonal
vectors $w_1,w_2$ of unit length in $\g\cdot v$ which span a $G_v$-section
in $\g\cdot v$; if the
$G$-action on $Q$ is polar, then there exists a section $\Si$ passing
through $\pi$ with $T_\pi\Si = {\text Span}(v \w w_1,v \w w_2)$. The
unique totally geodesic submanifold of $Q$ with such properties is given
by $\Si = \{v \w (ae + bw_1 + cw_2);\ a^2+b^2+c^2 = 1\}$. If we now
consider the two plane $\pi_1 = v \w w_1$, then $T_{\pi_1}\Si$ is
generated by $v \w e, v \w w_2$, while $T_{\pi_1}G\cdot\pi_1$ is given by
$\g v \w w_1 + v \w \g w_1$.
Since the scalar product $<v \w \g w_1,v \w w_2> = <\g w_1,w_2>$ is not
zero, we see that $\Si$ cannot be a section and this contradiction proves
our last claim.
\qed\enddemo
\bigskip
We now suppose that $G$ acts on $\R^n$ reducibly without fixed vectors and
that it acts on the quadric $Q_{n-2}$ coisotropically. We
split $\R^n = V_1 \oplus V_2$ into two nontrivial invariant subspaces
with $\dim V_i \geq 2$ and consider first the case when
both $V_i$ are $G$-irreducible. The restriction of the $\g$-action onto
$V_i$ 
defines homomorphisms $\rho_i: \g \to
\so(V_i)$ for $i=1,2$ .
We put $\g_i:= \rho_i(\g)$, $i=1,2$, so that $\g \subset \g_1 + \g_2$.
If $\dim V_i \geq 5$,   Proposition 4.3 and the fact that the action of $G$
restricted to
$V_i$ induces a coisotropic action on the corresponding quadric
imply that $\g_i$ is one of the following Lie algebras
$$\so(p), {\frak G}_2, \spin(7), \spin(9), \su(p), \u(p),
\R + \sp(p), \sp(1) + \sp(p) \tag *$$
for some integer $p$. On the other hand, if $2\leq \dim V_i \leq 4$, the
list in (*) contains all possible subalgebras of $\so(V_i)$ which
act on $V_i$ irreducibly, so that (*)
gives the list of all possible candidates for $\g_i$, $i=1,2$.
\medskip
\proclaim{Lemma 4.5} We have $\g_i \not= \su(p), \u(p), \R +
\sp(p)$.\endproclaim
\medskip
\demo{Proof} We suppose that $\g_1 \in \{\su(p),\u(p), \R + \sp(p)\}$. We
fix
 an oriented two plane $\pi$ in $V_2$ and let $\g_\pi$ denote
the stability subalgebra. If we identify $T_\pi Q$ with
$Hom(\pi,\pi^\perp)$,
we see that the slice representation at $\pi$ contains a submodule
isomorphic to $Hom(\pi,V_1)$. Let $J_0$ and $J_1$ denote $\g_\pi$-invariant
complex structures in $\pi$ and $V_1$ respectively and denote by $W_k$,
$k=1,2$,
the $\g_\pi$-submodules of $Hom(\pi,V_1)$ consisting of  elements $f$
such that $f(J_0v)=(-1)^kJ_1f(v)$ for all
$v\in V_1$. It is easy to see that $W_1$ and $W_2$ are equivalent
$\g_\pi$-modules. This implies that the $G$-action cannot be coisotropic
by Proposition 1.6.\qed\enddemo
\medskip
We now prove the following
\medskip
\proclaim{Lemma 4.6} If $\g$ is not isomorphic to $\g_1$, then $\g_2 \subset
\g$.
\endproclaim
\medskip
\demo{Proof} Since $\g$ is not isomorphic to $\g_1$, the kernel of $\rho_1$,
given by $\g\cap \g_2$, is not trivial. Therefore we can write $\g =
{\hat\g}_1
 + (\g\cap \g_2)$, sum of two ideals of $\g$, where ${\hat \g}_1$ is
isomorphic
 to $\g_1$. Now $\rho_2(\g\cap\g_2)$ is an ideal of $\g_2$; if $\g_2$ is
 simple, then $\rho_2(\g\cap\g_2) = \g\cap\g_2 = \g_2$, hence $\g_2\subseteq
\g$.
 If $\g_2$ is not simple, namely $\g_2\simeq \sp(1)+\sp(p)$ and $\g_2
 \not\subset\g$, then $\rho_2({\hat \g}_1)$ is a nontrivial ideal of
$\g_2$.
 We note that $\rho_2({\hat \g}_1)$ cannot coincide with $\g_2$: indeed, if
 so, we have $\g_1 \simeq \sp(1) + \sp(p) \simeq \g_2$ and $\rho_2|_{
 {\hat \g}_1}$ is an isomorphism; since $\ker \rho_2 = \g \cap \g_1$ is
 an ideal of $\g$ and $\g\cap \g_1\cap {\hat \g}_1 = \{0\}$, we have
 $\ker \rho_2 = \{0\}$ and $\dim \g = \dim \g_2 = \dim \g_1$, forcing
 $\g \simeq \g_1$, a contradiction. \par
 Therefore we have the following possibilities:\par
 \indent (a)\ $\rho_2({\hat \g}_1) = \sp(1)$. Then $\g\cap \g_2 = \sp(p)$.
 Now ${\hat \g}_1\simeq \g_1$ contains an ideal isomorphic to $\sp(1)$,
hence
 we have the following subcases:
\roster
\item "(a1)" $\g = \sp(1) + \sp(p)$, where $\rho_1(\g) = \sp(1)$ acting
on $\R^3$ by Lemma 4.5, $\rho_2(\g) = \sp(1) + \sp(p)$;
\item "(a2)" $\g = \sp(1) + \sp(p) + \sp(q)$, where $\g_1 = \sp(1)+\sp(q)$
and $\g_2 = \sp(1)+\sp(p)$.\endroster
\indent (b)\ $\rho_2({\hat \g}_1) = \sp(p)$. Then we have the following
subcases:
\roster
\item "(b1)" $\g = \sp(1)+\sp(p)$, where $\g_1 = \sp(p)$, $p=1,2$;
\item "(b2)" $\g = \sp(1)+\sp(p)+\sp(1)$ and $\g_1 \simeq \g_2$.
\endroster
We now show that the above cases do not yield coisotropic actions.
Indeed, for
(a1),
if we fix a two plane $\pi \subset \R^3=V_1$, then $\rho_2(\g_\pi) \subset
\R + \sp(p)$ leaves a complex structure invariant and the same arguments
as in Lemma 4.5 can be applied. Cases (a2) and (b2) can be dealt with
similarly.
As for (b1) a simple computation shows that the cohomogeneity $c$ of $G$
is at least $4$ if $p=1$ and $9$ if $p=2$, contradicting the fact
that $c$ does not exceed the rank of $G$.\qed\enddemo
\medskip
Therefore, if $\g$ is not isomorphic to $\g_i$
for $i=1,2$, then
$\g = \g_1 + \g_2$.
\medskip
\proclaim{Lemma 4.7} If $\g = \g_1 + \g_2$, then
\roster
\item "(i)" the action of $G$ on $Q_{n-2}$ is coisotropic if and only if
$$\g_i \in \{\so(V_i),{\frak G}_2,\spin(7),\spin(9),\sp(1) + \sp(p)\}$$
for $i = 1,2$;
\item "(ii)" the action of $G$ on $Q_{n-2}$ is polar if and only if
$$\g_i \in \{\so(V_i),{\frak G}_2,\spin(7)\}$$ for $i = 1,2$. \endroster
\endproclaim
\medskip
\demo{Proof} (i) We use the list (*) together with Lemma 4.5.\par Now if
$\g_i \in \{\so(p), {\frak G}_2, \spin(7)\}$ for
$i = 1,2$, then it can be easily checked
that $G$ has the same orbits in the quadric $Q = Q_{n-2}$ as
$\SO(V_1) \times \SO(V_2)$, which acts on $Q$ hyperpolarly, hence
coisotropically. \par
Now we suppose that, say, $\g_1 = \spin(9)$. We choose an oriented two plane
$\pi\subset V_2$ such that the orbit $G\cdot \pi$ is complex; then the
stability subalgebra $\g_\pi = \spin(9) + (\g_2)_\pi$ and
${\text{rank}} (\g_2)_\pi = {\text{rank}} (\g_2)$ since $G\cdot \pi$ is a
flag manifold. \par
We first claim that the $\g_\pi$-action on $\pi$
is not trivial.\par
Indeed, when $\g_2$ is $\so(V_2), {\frak G}_2$ or $\spin(7)$,
then our claim follows from the fact the groups $\SO(V_2), \G2$ and
$\Spin(7)$
act transitively on the Stiefel manifold of two frames in $V_2$ (see e.g.
[On], p.~90, 93). When $\g_2 = \spin(9)$, then the action of $\Spin(9)$ on
$Q_{14}$ has cohomogeneity one and we claim
that there exists only one complex singular
orbit in $Q_{14}$ given by $G\cdot\pi = \Spin(9)/\U(4)$, where $\U(4)$
acts nontrivially on $\pi$. Indeed, we see the spin representation of
$\Spin(9)$ as the isotropy representation of the Cayley plane
$\PO = \F4/\Spin(9)$
at the origin $o$ and, if we choose a two plane $\pi$ inside the tangent
space $W\cong \R^8$ at $o$ of a projective line in $\PO$ through $o$, then
the
stabilizer $\Spin(9)_\pi$ leaves $W$ invariant, since two different
projective lines cannot be tangent along a two plane; this means that
$\Spin(9)_\pi$ is contained in $\Spin(8)$ and therefore it covers
$\T^1\times \SO(6)\subset \SO(W)$. Hence $\Spin(9)_\pi = \U(4)$ and
it acts nontrivially on $\pi$. Moreover, since $\chi(Q_{14}) = 16$
is the sum of the Euler characteristics of the singular orbits in $Q_{14}$,
and since $\chi(\Spin(9)/\U(4)) = 16$, we see that there is
exactly one complex singular orbit $G\cdot\pi = \Spin(9)/\U(4)$.
When $\g_2 = \sp(1)+\sp(p)$, then the action of $\Sp(1)\cdot\Sp(p)$ on
$Q_{4p-2}$ has cohomogeneity one; if $\pi$ is a complex line, viewed as
a real two plane, inside a
quaternionic line in ${\Bbb H}^p$, then it is easy to see that
$(\g_2)_{\pi} = \R^2 + \sp(p-1)$ acts nontrivially on $\pi$.
The same argument as in
the previous case shows that $\Sp(1)\cdot\Sp(p)$ has only one
complex singular orbit in $Q_{4p-2}$.\par
So we may write $(\g_2)_\pi = \R + \h_2$, where $\h_2$ denotes the kernel of
the $(\g_2)_\pi$-action on $\pi$. The slice representation of $(\g_2)_\pi$
on $N_\pi$ splits as $(\R^2)^*\otimes V_1 + N'$, where $N'$ is acted on by
$\R + \h_2$; note that $N' = \{0\}$ if and only if $\g_2 \in
\{\so(p),{\frak G}_2, \spin(7)\}$. Now $(\R^2)^*\otimes V_1$ is acted on by
$\R + \spin(9)$, while $\h_2$ acts trivially on it; the representation of
$\R + \spin(9)$ on $(\R^2)^*\otimes V_1$ is coisotropic, see [Ka] or notice
that a regular
subalgebra is given by $\su(3)$ of corank three and the representation
 has cohomogeneity
three (see [Ya], p.~324).
If $\g_2 \in \{\so(p),{\frak G}_2, \spin(7)\}$, then
the $G$-action on the quadric $Q$ is coisotropic; indeed, a regular
isotropy subalgebra of $\g$ is $\su(3) + \h_2$ of corank three in
$\g_\pi$, hence of corank three in $\g$, while $N' = \{0\}$ implies that the
cohomogeneity
of the $G$-action on $Q$ is the cohomogeneity of $\R + \spin(9)$ acting on
$(\R^2)^*\otimes V_1$.\par
If $\g_2 = \spin(9)$, then $(\g_2)_\pi = \R + \su(4)$,
$N' = \R^8$ and $\su(4)$ acts on it by cohomogeneity one; hence a regular
isotropy subalgebra of $\g$ is $\su(3) + \su(3)$ of corank four, while the
cohomogeneity of the $G$-action is also four.\par
If $\g_2 = \sp(1) + \sp(p)$ ($p \geq 2$), then $(\g_2)_\pi = \R^2 +
\sp(p-1)$
and $N' = \R^{4p-4}$. So $(\g_\pi) = \su(3) + \R + \sp(p-2)$ of corank four
in $\g$ and the cohomogeneity of the $G$-action on $Q$ is also four.\par
Next we suppose that $\g_1 = \sp(1) + \sp(p)$. If we argue in the same way
as
in the case $\g_2=\spin(9)$, we see that we only need to prove
that the
representation of $\R + \sp(1) + \sp(p)$ on $(\R^2)^*\otimes \R^{4p}$ is
coisotropic (see also [Ka]). This representation has cohomogeneity three
(see
[Ya], p.~318); a regular isotropy subalgebra $\k$ contains a copy of
$\sp(p-2)$ as
an ideal, so that we can write $\k = \k' + \sp(p-2)$ and a simple dimension
count shows that $\k' = \R$. So the corank of $\k$ in $\R + \sp(1) + \sp(p)$
is also three and the representation is coisotropic. \par
(ii) Using the same arguments as in case (i), we only need to show that the
representation of $\R + \spin(9)$ on $\R^2 \otimes \R^{16}$ and of
$\R + \sp(1) + \sp(p)$ on $\R^2\otimes \R^{4p}$ are not polar. Actually
these representations are neither isotropy representations of irreducible
symmetric
spaces nor appear in the list of irreducible polar representations which are
not
isotropy representations of symmetric spaces (see the table from [EH] in the
proof of Lemma 2.7).\qed\enddemo
\medskip
Next we suppose that $\g$ is isomorphic to $\g_1$. \par
If $\g_1$ is simple, then $\rho_1,\rho_2$ are isomorphisms and $\g_1\simeq
\g_2$.
In particular, $\g_1 \in \{\so(V_1),{\frak G}_2,\spin(7),\spin(9)\}$. Now,
if $\rho_1$ is equivalent to $\rho_2$, we select a two plane $\pi\subset
V_1$
and observe that the slice representation of $\g_\pi$ contains
$Hom(\pi,V_2)$
as a complex submodule; if we denote by $\pi$ the same oriented two plane
in $V_2\simeq V_1$, then $Hom(\pi,\pi)$ contains a trivial
$\g_\pi$-submodule.
This means that the slice representation contains a trivial complex
submodule,
preventing the $G$-action from being coisotropic. If on the other
hand $\rho_1$ is not equivalent to $\rho_2$, then $\g=\spin(7)$ or
$\g=\spin(9)$; in this case a simple computation shows that the dimension
of a Borel subalgebra of $\g^\C$ is strictly less than the complex
dimension of the quadric and the action is therefore not
coisotropic.\par
If $\g_1$ is not simple, i.e. $\g_1 = \sp(1)+\sp(p)$, then $\g_2$ can be
$\sp(1)\simeq\so(3)$, $\sp(2)\simeq \so(5)$ or $\sp(1)+\sp(p)$. The first
two cases already appeared in the proof of Lemma 4.6  and the last case
 can be ruled out with the same arguments.\enddemo\par

In order to complete the proof of Theorem 1.3, we need to prove the
following
lemma.
\medskip
\proclaim{Lemma 4.8} Let $G$ be a compact Lie group acting reducibly on
$\R^n$ ($n\geq 5$) and coisotropically on the quadric $Q_{n-2}$. Then the
number
of irreducible submodules of $\R^n$ is at most two.\endproclaim
\medskip
\demo{Proof} We assume that $\R^n$ contains at least three irreducible
submodules $V_1,V_2,V_3$, none of which is trivial. We then
restrict ourselves to the coisotropic action of $G$ on the quadric
corresponding to the $G$-module $V = \sum_{i=1}^3 V_i$. We can apply the
results
obtained so far to the quadrics corresponding to the sums $V_i+V_j$,
$i,j \in\{1,2,3\}$, showing that $\g = \g_1+\g_2+\g_3$, where
each $\g_i \in
\{\so(V_i),{\frak G}_2,\spin(7), \spin(9), \sp(1)+\sp(p)\}$. We consider a
two plane
$\pi\subset V_1$
and its isotropy subalgebra $\g_\pi = (\g_1)_\pi+\g_2+\g_3$. The
slice representation of $\g_\pi$ contains $Hom(\pi,V_2+V_3)$ as a complex
submodule which is coisotropic; this means that $\g_\pi$ acts nontrivially
on $\pi$, since otherwise $Hom(\pi,V_2+V_3)$ would contain equivalent
$\g_\pi$-submodules, contradicting Proposition 1.6.  The action of $\g_\pi$
on
$Hom(\pi,V_2+V_3)$ is the action of $\h:= \R+\g_2+\g_3$ on
$(\R^2)^*\otimes (V_2+V_3)$.
We set $W_i:= (\R^2)^*\otimes V_i$ and $\k_i$ the regular isotropy
subalgebras of $\R+\g_i$ acting on $W_i$ for $i = 2,3$; it is known that
the regular isotropy subalgebra of $\R + \so(n)$ acting on
$\R^2 \otimes \R^n$ is $\so(n-2)$, hence $\k_i\subset \g_i$ for $i = 2,3$.
Now a regular isotropy subalgebra $\k$ of $\h$ acting on
$Hom(\pi,V_2+V_3)$ is
given by $\k_1 + \k_2$. Since the $\R + \g_i$-actions on $W_i$ are
coisotropic, we have that the cohomogeneity $c$ of the $\h$-action
on $Hom(\pi,V_2+V_3)$ is given by
$$\eqalign{ c = &\;\dim (W_1 + W_2) - \dim(\g_1 + \g_2) - 1 +
\dim(\k_1+\k_2)\cr
{} = &\;\dim W_1 - \dim (\R + \g_1) + \dim \k_1 \cr
{}& \; + \dim W_2- \dim (\R + \g_2) + \dim \k_2 + 1 \cr
{} = &\;\; {\text{rk}}(\R + \g_1) - {\text{rk}} (\k_1) + {\text{rk}}(\R +
\g_2)
- {\text{rk}}(\k_2) + 1\cr
{} = &\;\;{\text{rk}}(\R + \g_1 + \g_2) - {\text{rk}}(\k_1 + \k_2) + 2.}$$
This means that the $G_\pi$ action on $Hom(\pi,V_2+V_3)$ is not coisotropic.
Hence the
$G$-action on $Q_{n-2}$ is not coisotropic. \qed\enddemo

\bigskip

\hbox{\parindent=0pt\parskip=0pt
\vbox{\hsize=2.7truein
\obeylines
{
Fabio Podest\`a
Dip.~di Matematica e Appl.~per l'Arch.
Universit\`a di Firenze
P.zza Ghiberti, 27
I-50142 Firenze
Italy
}\medskip
podesta\@math.unifi.it
}\hskip 1.5truecm
\vbox{\hsize=3.7truein
\obeylines
{
Gudlaugur Thorbergsson
Mathematisches Institut
Universit\"at zu K\"oln
Weyertal 86-90
D-50931 K\"oln
Germany
}\medskip
gthorbergsson\@mi.uni-koeln.de
}
}

\enddocument
\bye